\documentclass[12pt,reqno]{amsart}
\usepackage{amsmath}
\usepackage{amscd,amsthm,amsfonts,amsopn,amssymb,mathrsfs, latexsym}
\usepackage{epsfig,hyperref}
\theoremstyle{plain}
\newtheorem{theorem}{Theorem}

\newtheorem{proposition}[theorem]{Proposition}

\theoremstyle{remark}

\allowdisplaybreaks
\def\be{\begin{equation}}
\def\ee{\end{equation}}
\parskip=\medskipamount

\def\vp{\varphi}

\def\arrowk{^\to{\kern -6pt\topsmash k}}
\def\arrowK{^{^\to}{\kern -9pt\topsmash K}}
\def\arrowr{^\to{\kern-6pt\topsmash r}}
\def\arrowvp{^\to{\kern -8pt\topsmash\vp}}
\def\arrowf{^{^\to}{\kern -8pt f}}
\def\arrowg{^{^\to}{\kern -8pt g}}
\def\arrowu{^{^\to}a{\kern-8pt u}}
\def\arrowt{^{^\to}{\kern -6pt t}}
\def\arrowe{^{^\to}{\kern -6pt e}}
\def\tk{\tilde{\kern 1 pt\topsmash k}}
\def\barm{\bar{\kern-.2pt\bar m}}
\def\barN{\bar{\kern-1pt\bar N}}
\def\barA{\, \bar{\kern-3pt \bar A}}

\def\iint{\not \kern-4pt\int}
\begin{document}
\title{PARTIAL QUOTIENTS AND REPRESENTATION OF RATIONAL NUMBERS\\
\ \\ QUOTIENTS PARTIELS ET REPR\'ESENTATION DES NOMBRES RATIONNELS}

\author{Jean Bourgain}
\address{School of Mathematics, Institute for Advanced Study, 1
Einstein Drive, Princeton, NJ 08540.}
\email{bourgain@math.ias.edu}
\thanks{The research was partially supported by NSF grants DMS-0808042 and
DMS-0835373.}
\maketitle

\begin{abstract}
It is shown that there is an absolute constant $C$ such that any rational $\frac bq\in ]0, 1[, (b, q)=1$, admits a representation as a finite sum $\frac bq =\sum_\alpha \frac {b_\alpha}
{q_\alpha}$ where $\sum_\alpha \sum_i a_i \big(\frac {b_\alpha}{q_\alpha}\big)< C \log  q$ and $\{a_i(x)\}$ denotes the sequence of partial quotients of $x$.
\end{abstract}

\renewcommand{\abstractname}{R\'esum\'e}
\begin{abstract}
On d\'emontre l'existence d'une constante $C$ telle que tout rationnel $\frac bq\in ]0, 1[, (b, q)=1$, a une repr\'esentation comme somme finie $\frac bq= \sum_\alpha \frac {b_\alpha}{q_\alpha}$
o\`u $\sum_\alpha\sum_i a_i\big(\frac {b_\alpha}{q_\alpha}\big)< C\log q$ et $\{a_i(x)\}$ est la suite des quotients partiels de $x$.
\end{abstract}

\bigskip
\section*{Version fran\c caise abr\'eg\'ee}

Cette Note r\'esulte de r\'esultats r\'ecents sur la conjecture de Zaremba obtenus dans \cite {B-K} et quelques questions pos\'ees par R.~Kenyon \cite {K} sur la repr\'esentation de nombres
rationnels comme somme de nombres rationnels dont les quotients partiels sont born\'ees par une constant absolue.
Diverses propri\'et\'es de repr\'esentation de nombres r\'eels $x$ comme somme $x=y+z+\cdots$, o\`u $y, z, \ldots$ ont leurs quotients partiels sujets \`a certaines bomes, ont en effet \'et\'e
\'etablies (voir en particulier les r\'esultats de M.~Hall, \cite {H}).
Le probl\`eme de trouver des \'enonc\'es analogues pour les rationnels semble donc naturel.
Dans cet esprit, on exploite ici les m\'ethodes de \cite {B-K} afin d'\'etablir la propri\'et\'e suivante.

\begin{proposition}
Il existe une constant $C$ telle que tout rationnel \hfill\break
$\frac bq \in ]0, 1[, (b, q)=1$, admette une repr\'esentation comme somme finie $\frac bq =\sum_\alpha\frac {b_\alpha}{q_\alpha}$
o\`u $\sum_\alpha \sum_i a_i \big(\frac {b_\alpha}{q_\alpha}\big)<C\log q$.
\end{proposition}

\section 
{Some background}

It was shown by M.~Hall \cite{H} that every number in the interval \hfill\break
$]\sqrt 2-1, 4\sqrt 2-4[$ is the sum of two continued fractions whose partial quotients do not exceed four (see \cite{H}, Th.
3.1).

Recently, R.~Kenyon brought to the author's attention the problem of obtaining a result in the flavor of Hall's theorem for the rational numbers.
There are several possible formulations.
One could ask for instance if there is an absolute constant $C$ such that given $\frac bq\in\mathbb Q_+\cap I$, $I$ a suitable interval, there is a representation of $b$ as a sum of at most
$C$ positive integers $b_i$ such that each of the fractions $\frac {b_i}q$ has its partial quotients bounded by $C$.
To be mentioned here is the (still unsolved) conjecture of Zaremba, according to which for all $q\in\mathbb Z_+$, there is some $(b, q)=1$ such that $\frac bq$ has partial quotients bounded
by five (or some absolute constant).
Alternatively, one can ask if any element in $\mathbb Q_+\cap I$ is sum of two (or at most $C$) rationals with partial quotients bounded by $C$.
While we will leave these questions unanswered here, our aim is to prove the following property in a similar spirit.

\noindent
{\bf Proposition 1.}
{\it There is an absolute constant $C$ such that any rational $\frac bq\in]0, 1[, (b, q)=1$, admits a representation as a finite sum
\begin{equation}\label{1.1}
\frac bq =\sum_\alpha \pm \frac {b_\alpha}{q_\alpha} \qquad (b_\alpha, q_\alpha)=1
\end{equation}
such that
\begin{equation}\label{1.2}
\sum_\alpha\sum_i a_i\Big(\frac {b_\alpha}{q_\alpha}\Big) \leq C\log q
\end{equation}
where $\{a_i(x)\}$ denotes the sequence of partial quotients of $x\in ]0, 1[$.}

As R.~Kenyon points out, a statement of this kind may be viewed as a measure of complexity of rationals of given height.
Note that since $\sum_i a_i\big(\frac {b_\alpha}{q_\alpha}\big)\gtrsim \log q_\alpha$, above estimate is essentially optimal.

\section
{Preliminaries}

Our main analytical tools are the results and methods of the recent paper \cite{B-K} on Zaremba's conjecture.
It is shown in \cite{B-K} that for a large enough constant $A$ (we may take $A=50$), for all $q\in\mathbb Z_+$ outside an exceptional set $E\subset\mathbb Z_+$ of zero-density, there is
some $b\in\mathbb Z_+$, $(b, q)=1$ such that
\begin{equation}\label{2.1}
\frac bq\in\mathcal R_A =\{x\in\mathbb Q\cap [0, 1]; \max_i a_i(x)\leq A\}.
\end{equation}
More quantitatively, one gets an estimate
\begin{equation}\label{2.2}
|E\cap [1, N]|< N^{1-\frac c{\log\log N}}.
\end{equation}
Instead of \eqref{2.2}, it is possible to obtain a power saving

\noindent
{\bf Proposition 2.}
\label{proposition2}
{\it The above statement holds with $E$ satisfying
\begin{equation}\label {2.3}
|E\cap [1, N]|< N^{1-c_1}
\end{equation}
with $c_1>0$ some constant.}

Recalling the approach from \cite{B-K}, elements $\frac bq\in\mathcal R_A$ are produced from elements $g=\begin{pmatrix} *&b\\ *&q \end{pmatrix}$ in the semi-group $\mathcal G_A$ generated by the
matrices
\begin{equation}\label {2.4}
\begin{pmatrix} 0& 1\\ 1&a\end{pmatrix} \qquad (1\leq a\leq A).
\end{equation}
We use the Hardy-Littlewood circle method in order to analyze exponential sums of the form
\begin{equation}\label{2.5}
\sum \lambda(g) 
e(g_{22} \theta)
\end{equation}
with $\lambda$ a suitable distribution on $\mathcal G_A$.
Fixing some large $N$, the distribution $\lambda$ is obtained from a product of certain Archimedian
balls in $\mathcal G_A$.
As usual, the circle method involves a treatment of minor and major arcs contributions and those depend on different ingredients.
The estimates on minor arcs result from Vinogradov-type bilinear estimates, exploiting the multi-linear structure of $\lambda$.
A precise evaluation of \eqref{2.5} on the major arcs (up to an error term) is possible using spectral methods.
We use the thermodynamical approach and the results from \cite{BGS} based on the theory of expansion in $SL_2(q)$.
The error term in the counting and the size of the exceptional set in \eqref{2.2} depend on the width of the resonance free 
regions for the congruence transfer operators.
It turns out that the gain of the $N^{-\frac c{\log\log N}}$-factor rather than $N^{-c}$ comes from introducing balls $B_M=\{g\in\mathcal G_A; \Vert g\Vert \leq M\}$.
If instead of balls we consider slightly more general distributions (obtained as average of balls over suitable radii), one may recover a full powergain $N^{-c}$ (the cutoff-level for the major
arcs may then be set at $N^c$ for some $c>0$, rather than $N^{\frac c{\log\log N}}$).

In order to prove Proposition 1, we need one more further refinement.

\noindent
{\bf Proposition 3.}
{\it Taking again $A$ sufficiently large, there is a constant $c>0$ such that the following holds.

Let $N\in\mathbb Z_+$ be large enough, $d\in\mathbb Z_+, d<N^c$ and  $\beta\in\mathbb Z$, $(\beta, d)=1$.
There is a subset $E_{N; d, \beta}\subset\mathbb Z\cap [1, N]$ such that
\begin{equation}\label{2.6}
|E_{N, d, \beta}|<N^{1-c}
\end{equation}
and for all $q\in\mathbb Z_+\backslash E_{N; d, \beta}, \, q<N$, there is $b\in\mathbb Z_+, b< q, (b, q)=1$ satisfying
\begin{equation}\label{2.7}
\frac bq\in\mathcal R_A
\end{equation}
and
\begin{equation}\label{2.8}
b\equiv \beta \ (\text{mod\,} d).
\end{equation}
}

Returning  to \cite{B-K}, the incorporation of the additional congruence condition \eqref{2.8} is harmless at the level of the minor arcs estimates,
provided $d$ is sufficiently small.
Of course the condition (10) enters the singular series in the treatment of the major arcs and the assumption $(\beta, d)=1$ ensures that there are no local obstructions.

Obviously Proposition 3 implies that for some constant $c>0$, the following holds

\noindent
{\bf Proposition 3$'$.}
{\it There is a subset $E_N\subset\mathbb Z\cap [1, N]$ such that
\begin{equation}\label{2.9}
|E_N|< N^{1-c}
\end{equation}
and for all $q\in\mathbb Z_+\backslash E_N, q< N$ and all $d\in\mathbb Z_+, d<N^c$, and $\beta\in\mathbb Z$,  $(\beta, d)=1$, there is some $b\in\mathbb Z_+, (b, q)=1$ satisfying \eqref {2.7} and \eqref{2.8}.
}

\section
{Proof of Proposition 1}

Denote $C(\frac bq)$ the minimum of the left hand side of \eqref{1.2} over all representations \eqref{1.1}.
First, observe that it suffices to show that
\begin{equation}\label{3.1}
\frac bq=\frac {b'}{q'}+\frac {b''}{q''}
\end{equation}
with
\begin{equation}\label{3.2}
C\Big(\frac {b'}{q'}\Big)<C\log q
\end{equation}
and
\begin{equation}\label{3.3}
q'' <\sqrt q
\end{equation}
with $C$ in \eqref{3.2} some absolute constant.
We may then indeed iterate.

Let $c>0$ be the constant from Proposition 3$'$. Set
\begin{equation}\label {3.4}
\delta=\frac 1{10} c \text { and } r=[\delta^{-2}]+1.
\end{equation}
We claim that there are primes $p_1, \ldots, p_r$, $(p_i, q)=1$, $p_i< q^\delta$ and $p_i\sim q^\delta$ satisfying the following two
conditions

\begin{equation}\label {3.5}
qp_1\ldots p_r \not\in E_{q^{1+r\delta}}\text { with $E_N$ as in Proposition 3$'$}
\end{equation}

\begin{equation}\label{3.6}
\text {For all $I\subset\{1, \ldots, r\}, I\not= \phi, \prod_{i\in I} p_i\not\in E_{q^{\delta|I|}}$}.
\end{equation}
Indeed, consider all integers of the form $qp_1\ldots p_r$ with $p_i$ as above.

Their number is at least
\be\label{3.7}
\frac {q^{r\delta}}{(\log q)^r}.
\ee
On the other hand, by \eqref{2.9}
$$
|E_{q^{1+r\delta}}|< q^{(1+r\delta)(1-c)} \leq q^{1+r\delta-c\delta^{-1}} < q^{r\delta-1}=o\Big(\frac{q^{r\delta}}{(\log q)^r}\Big).
$$
Next, consider condition \eqref{3.6} and fix $I\subset\{1, \ldots, r\}, I\not=\phi$.
Among the integers considered above, those for which $\prod_{i\in I} p_i\in E_{q^{|I|\delta}}$ account for at most
$$
q^{(r-|I|)\delta} |E_{q^{\delta|I|}}|  < q^{r\delta-c\delta} =o\Big(\frac {q^{r\delta}}{(\log q)^r}\Big).
$$
Hence we may find $p_1, \ldots, p_r$ with the desired properties.

Returning to \eqref{3.1}-\eqref{3.3}, write with $p_1, \ldots, p_r$ as above
$$
\frac bq=\frac {bp_1\ldots  p_r} {qp_1\ldots p_r}.
$$
Since $(b, q) =(p_i, q)=1$, $(bp_1\ldots p_r, q)=1$ with $q< N^c$, $N=q^{1+r\delta}$.

Since \eqref{3.5} holds, Proposition 3$'$ implies that there is some $b_0\in\mathbb Z_+$, $b_0< qp_1\ldots p_r$ such that $\frac{b_0}{p_1\ldots p_rq}\in\mathcal R_A$ and
$b_0\equiv bp_1\ldots p_r (\text{mod\,} q)$.
Hence
\be\label{3.8}
C\Big(\frac {b_0}{p_1\ldots p_r q}\Big) < C(A) \log (p_1\ldots p_rq)\leq C(A) (1+r\delta) \log q =C\log q
\ee
and we may write
\be\label{3.9}
\frac bq =\frac {b_0}{qp_1\ldots p_r} +\frac {a_1}{\prod_{i\in I_1} p_i} 
\text { with } a_1\in\mathbb Z  \text { and } \Big( a_1, \prod_{i\in I_1}p_i\Big) =1 \text { if } I_1\not=\phi.
\ee
If $\prod_{i\in I_1}p_i<\sqrt q$, set $\frac {b''}{q''} =\frac {a_1}{\prod_{i\in I_1} p_i}$ in \eqref{3.1}.

If $\prod_{i\in I_1} p_i\geq \sqrt q$, use \eqref{3.6} and take $d=p_{i_1}$, $i_1$ chosen from $I_1$, $\beta=a_1$.

By definition of $\delta, p_{i_1} < q^\delta < (q^{\frac 12})^c$ and we get some $b_1\equiv a_1(\text{mod\,} p_{i_1})$ with $\frac {b_1}{\prod_{i\in I_1}
p_i}\in \mathcal R_A$. Thus
\be\label{3.10}
C\Big(\frac {b_1}{\prod_{i\in I_1}p_i}\Big) < C(A) \log\Big(\prod_{i\in I_1}p_i\Big) < C\log q
\ee
and

\be\label{3.11}
\frac {a_1}{\prod_{i\in I_1} p_i} =\frac {b_1}{\prod_{i\in I_1}p_i} +\frac {a_2}{\prod_{i\in I_2}p_i} \text { with } I_2\subset I_1\backslash \{i_1\},
\Big(a_2,\prod_{i\in I_2}p_i\Big)
=1 \text { if } I_2\not=\phi.
\ee

The continuation of the process is clear and it terminates after at most $r$ steps, leading to a representation
\be\label{3.12}
\frac bq=\frac {b_0}{qp_1\ldots p_r}+\frac {b_1}{\prod_{i\in I_1}p_i} +\cdots+\frac {b_\rho}{\prod_{i\in I_\rho} p_i} +\frac {b''}{q''} =\frac {b'}{q'}+\frac {b''}{q''}
\ee
satisfying \eqref{3.2}, \eqref{3.3}.
\bigskip

\noindent
{\bf Acknowledgement.}  The author is grateful to R.~Kenyon for bringing these questions to his attention.

\end{document}